\documentclass{hha}
\usepackage{euscript}

\newtheorem{theorem}{Theorem}[section]
\newtheorem{definition}[theorem]{Definition}
\newtheorem{problem}[theorem]{Problem}
\newtheorem{observation}[theorem]{Observation}
\newtheorem{example}[theorem]{Example}
\newtheorem{lemma}[theorem]{Lemma}

\newtheorem{proposition}[theorem]{Proposition}
\newtheorem{remark}[theorem]{Remark}

\def\pa{\partial} \def\phi{\varphi}
\def\catdgVect{{\tt dgVect}}

\def\am{{A($m$)}}\def\an{{A($n$)}}
\def\lm{{L($m$)}}
\def\catam{{\tt A}($m$)} \def\catan{{\tt A}($n$)}
\def\catlm{{\tt L}($m$)} \def\catln{{\tt L}($n$)}

\def\bk{{\bf k}} 

\def\rada#1#2#3{#1_{#2},\dots,#1_{#3}}

\def\sumdeg#1#2#3{|#1_{#2}|+\dots+|#1_{#3}|}
\def\ucalam{{\hskip .2em\raisebox{.72em}{\tiny $m$} \hskip -.8em 
         \underline{\EuScript A}}}
\def\ucalamdg{{\hskip .2em\raisebox{.72em}{\tiny $m$}\hskip -.8em 
          \underline{\mathcal A}}}

 \def\varomega{\vartheta}
\def\uGamma{\underline{\Gamma}}
\def\ucalfm{{\hskip .2em\raisebox{.72em}{\tiny $m$}\hskip-.72em 
          \underline{\EuScript F}}}
\def\calfm{{\hskip .2em\raisebox{.72em}{\tiny $m$}\hskip-.72em {\EuScript F}}}
\def\callm{{\hskip .2em\raisebox{.72em}{\tiny $m$} \hskip -.7em {\EuScript L}}}
\def\callff#1#2{{\hskip .2em \raisebox{.72em}{\scriptsize 
          $#1 \hskip -.8mm :  \hskip -.6mm#2$} \hskip -.9em {\EuScript L}}}
\def\calldve{{\buildrel \mbox{\scriptsize $2$} \over {\EuScript L}}}

\def\treem{{\raisebox{.72em}{\tiny $m$}\hskip-.58em {\mbox{$\tt T$}}}}
\def\Span{{\rm Span}_{{\bf k}}}
\def\calim{{\hskip .2em\raisebox{.72em}{\tiny $m$}\hskip -.6em {\EuScript I}}}
\def\ucalim{{\hskip .2em\raisebox{.72em}{\tiny $m$}\hskip -.6em 
          \underline{\EuScript I}}}

\def\frm{{\raisebox{.72em}{\tiny $m$}\hskip -.7em {\mathfrak A}}}
\def\frmff#1#2{{\raisebox{.72em}{\tiny $#1 \hskip -.8mm :  \hskip -.6mm#2$} 
         \hskip -.9em{\mathfrak A} \hskip .2em}}
\def\frmfff#1{{\hskip .2em\raisebox{.72em}{\tiny  $#1$} 
               \hskip -.5em {\mathfrak A}}}
\def\frml{{\hskip .2em\raisebox{.72em}{\tiny $m$} \hskip -.7em  {\mathfrak L}}}
\def\vert{{\rm Vert}}
\def\val{{\rm val}}
\def\pphi#1{{\hskip .3em\raisebox{.72em}{\tiny $#1$} \hskip -.7em {\phi}}}
\def\vvaromega#1{{\buildrel \mbox{\scriptsize $#1$} \over {\varomega}}}
\def\VECTOR#1{\vec #1}
\def\bm{{{\hskip .2em\raisebox{.72em}{\tiny $m$} \hskip -.8em B}}}
\def\adm#1{{\hskip .8em\raisebox{.72em}{\tiny $#1$}\hskip -1.4em{{\rm Adm}}}}
\def\admm{\adm m}
\def\b{\bullet}
\def\c{\circ}
\def\LLL#1{\mbox{\rm L($#1$)}}
\def\sgn{{\rm sgn}}
\def\ttt{{\tt T}}
\def\ul{{\mathcal U}(L)}
\def\ot{{\otimes}}
\def\calu{{\mathcal U}}

\def\AM{{\mbox{{\rm A(}$m${\rm )}}}}

\def\cals{{\EuScript S}}

\def\vert{{\rm Vert}}

\begin{document}

\bibliographystyle{plain}

\title{Free homotopy algebras}

\author{Martin Markl}
\email{markl@math.cas.cz}
\address{Mathematical Institute of the Academy\\
         \v Zitn\'a 25\\
         115 67 Praha 1\\ The Czech~Republic}
\thanks{This work was supported by the
grant GA~\v CR 201/96/0310 and M\v SMT ME603.}


\begin{abstract}
  An explicit description of free strongly homotopy associative and
  free strongly homotopy Lie algebras is given. A variant of the
  Poincar\'e-Birkhoff-Witt theorem for the universal enveloping
  \am-algebra of a strongly homotopy Lie algebra is formulated.
\end{abstract}



\maketitle


\section*{Introduction}

This note was originated many years ago as my reaction to questions of
several people how free strongly homotopy algebras can be described
and what can be said about the structure of the universal enveloping
\am-algebra of an \lm-algebra constructed
in~\cite{lada-markl:CommAlg95}, and then circulated as a ``personal
communication.'' I~must honestly admit that it contains no really deep
result and that everything I did was that I expanded definitions and
formulated a couple of statements with more or less obvious proofs.

\am-algebras and their strict
homomorphisms~\cite[pages~147--148]{markl:JPAA92} form an equationally
given algebraic category ${\tt A}(m)$.  It follows from general theory
that the forgetful functor to the category ${\tt gVect}$ of graded
vector spaces, $\Box :{\tt A}(m) \to {\tt gVect}$, has a left adjoint
$\frm : {\tt gVect}\to {\tt A}(m)$. Given a graded vector space $X \in
{\tt gVect}$, the object $\frm(X)\in {\tt A}(m)$ is the {\em free}
\am{\em -algebra\/} on the graded vector space $X$. We will also, for
$n < m$, consider forgetful functors $\Box : {\tt A}(m) \to{\tt A}(n)$
and their left adjoints $\frmff mn : {\tt A}(n) \to{\tt A}(m)$; here
the case $n=1$ is particularly important, because $\frmff m1 :
\catdgVect \to \mbox {\catam}$ describes the free \am-algebra
generated by a {\em differential\/} graded vector space.

We believe there is no need to emphasize the important r\^ole 
of free objects in mathematics. Each \am-algebra is a quotient of a
free one and free \am-algebras were used in our definition of the
universal enveloping algebra of a strongly homotopy 
Lie algebra~\cite[page~2154]{lada-markl:CommAlg95}. It
could be useful to have a concrete description of these free
algebras, as explicit as, for example,
the description of free associative algebras by tensor algebras.

This brief note gives such a description in terms of planar trees. 
Replacing planar trees by
non-planar ones, one can equally easy represent also {\em free 
strongly homotopy Lie algebras}.
Surprisingly, these free algebras are
simpler objects that their strict counterparts
and admit a nice linear basis 
(Section~\ref{to_jsem_zvedav_jestli_aspon_neco_bude_prijato}), while 
free (strict) Lie algebras are very complicated 
objects
(see, for example,~\cite[Chapter~IV]{serre:65}).
This might be explained by the fact that axioms of strongly 
homotopy algebras are certain resolutions of
strict axioms, thus they are more amenable. They
 manifest some properties of distributive
laws~\cite{markl:dl} -- they are `directed,' therefore one can
obtain unique representatives of elements of the corresponding
operads in terms of admissible trees, like 
in~Proposition~\ref{bez_odezvy}, with obvious implications for the
structure of free algebras.

In the first section of this note we describe the operad $\ucalam$ for
\am-algebras\footnote{\underline{underlinin}g indicates that it is a
  non-$\Sigma$ operad} and then, in Section~2, free \am-algebras. In
the third section we modify these constructions to free strongly
homotopy Lie algebras. In the last section we prove a
Poincar\'e-Birkhoff-Witt-type theorem for strongly homotopy Lie
algebras and formulate a problem related to the structure of the
category of modules. In the appendix we then recall the language of
trees used throughout the paper.

\section{Operad for \am-algebras}

In this note we assume a certain familiarity with operads, a good
reference book for this subject is~\cite{markl-shnider-stasheff:book}.
For simplicity we suppose the ground ring $\bk$ to be a field of
characteristic zero, though the results of this and the following
section remain true over an arbitrary ring.

Let $m$ be a natural number or $\infty$. Recall 
(see~\cite{stasheff:TAMS63} for the original definition, 
but we use the sign convention of~\cite{markl:JPAA92})
that an \am-{\em algebra}
is a graded $\bk$-module $A$ together with a set $\{ \mu_k; 1\leq k
\leq m,\ k<\infty\}$ 
of degree $k-2$-linear maps, $\mu_k : \bigotimes^k A \to A$, such that
\[
\sum_{\lambda =0}^{n-1} \sum_{k=1}^{n-\lambda} 
(-1)^\omega \cdot
\mu_{n-k+1}(\rada a1\lambda,\mu_k(\rada
a{\lambda+1}{\lambda+k}),
\rada a{\lambda+k+1}n) = 0
\]
for all homogeneous $a_1,\ldots,a_n \in A$, $n \leq m$, where the sign
is given as
\[
\omega := {k+\lambda+k\lambda +k(\sumdeg a1{\lambda})}.
\]

It easily follows from the above axiom that $\mu_1 : A \to A$ is a
degree $-1$ differential. This operation will play a special r\^ole
and we denote it by $\partial := \mu_1$.  \am-algebras are algebras
over a certain non-$\Sigma$ operad $\ucalam= \{\ucalam(n) \}_{n\geq
  1}$ in the monoidal category of graded vector spaces.  This operad
can be constructed as follows.

Let $\ucalfm = \{\ucalfm(n)\}_{n\geq 1}$ be the free graded non-$\Sigma$
operad $\uGamma(\xi_1,\xi_2,\ldots,\xi_m)$ generated by operations
$\xi_i$, $l\leq i \leq m$, where each $\xi_i$ has arity $i$ and degree
$i-2$.  Alternatively, let $\treem(n)$ be the set of all connected
directed planar trees with $n$ input edges, all of whose vertices have
at most $m$ input edges; we admit also `unary' vertices with only one
input edge.

Let $\vert(T)$ denote the set of vertices of
a tree $T \in \treem(n)$. For $v\in \vert(T)$, let $\val(v)$ be the number of
input edges of $v$. We assign to each $T\in \treem(n)$ a degree
$\deg(T)$ by the formula
\[
\deg(T) := \sum_{v\in \vert(T)}(\val(v)-2).
\]
Then $\ucalfm(n)$ can be naturally identified with the $\bk$-linear
span of the set $\treem(n)$,
\begin{equation}
\label{1}
\ucalfm(n) = \Span(\treem(n)),
\end{equation}
see~\cite[Section~II.1.9]{markl-shnider-stasheff:book}.  The operadic
structure is, in this language, defined by grafting the corresponding
trees.  This means that, for $T\in \treem(a)$ and $S\in \treem(b)$,
the operadic $i$-th circle
product~\cite[Definition~II.1.16]{markl-shnider-stasheff:book} $T
\circ_i S$ is the tree obtained by grafting the tree $S$ at the $i$-th
input edge of $T$. Similarly, for $T\in \treem(l)$ and $S_j \in
\treem(k_j)$, $1\leq j \leq l$, the operadic
composition~\cite[Definition~II.1.4]{markl-shnider-stasheff:book}
$\gamma(T;S_1,\ldots,S_l) =T(S_1,\ldots,S_l) \in
\treem(k_1+\cdots+k_l)$ is the tree obtained by grafting $S_j$ at the
$j$-th input edge of $T$, for each $1\leq j \leq l$.

The operad $\ucalam$ is then defined as $\ucalfm/\ucalim$, where $\ucalim$
is the operadic ideal generated by the elements
\begin{equation}
\label{2}
\Phi_n := \sum_{\lambda=0}^{n-1}
\sum_{k=1}^{n-\lambda}(-1)^{\lambda(k+1)+k}
\xi_{n-k+1}\circ_{\lambda+1} \xi_k \in \ucalfm(n),\ n \leq m.
\end{equation}
In the tree language  of~(\ref{1}), the generators $\xi_j$
are represented by $j$-corollas (denoted by the same symbol)
$\xi_j \in \treem(j)$; recall that the $j$-corolla is the unique
tree having exactly one vertex. Relations~(\ref{2}) then represent
certain linear combinations of trees, indicated
in Figure~\ref{picture}.
\begin{figure}
\begin{center}
\unitlength 1.50mm
\thicklines
\begin{picture}(44.50,20.00)(0,17)
\put(22.00,25.00){\makebox(0,0)[cc]{$\Phi_n := \displaystyle%
\sum_{\lambda=0}^{n-1}\sum_{k=1}^{n-\lambda}(-1)^{\lambda(k+1) + k}%
\left(\rule{0mm}{17mm}\hskip39mm \right)$}}
\put(-5,0){\put(35.17,35.00){\line(0,-1){5.00}}
\put(35.17,30.00){\line(-6,-5){6.00}}
\put(35.17,30.00){\line(-1,-2){2.50}}
\put(35.17,29.67){\line(1,-3){3.11}}
\put(38.28,20.33){\line(-1,-1){6.11}}
\put(38.33,20.00){\line(-2,-5){2.27}}
\put(38.33,20.00){\line(1,-1){5.50}}
\put(35.17,29.83){\line(2,-1){9.33}}
\put(39.17,15.83){\makebox(0,0)[cc]{$\cdots$}}
\put(39.50,25){\makebox(0,0)[cc]{$\cdots$}}
\put(35,25){\makebox(0,0)[cc]{$\cdots$}}
\put(38.17,22.00){\makebox(0,0)[lc]{$\lambda+1$-th input}}
\put(35.17,29.8){\makebox(0,0)[cc]{$\bullet$}}
\put(38.33,20.00){\makebox(0,0)[cc]{$\bullet$}}
\put(42.17,32.50){\makebox(0,0)[cc]{$\xi_{n-k+1}$}}
\put(44.33,17.83){\makebox(0,0)[cc]{$\xi_k$}}}
\end{picture}
\end{center}
\caption{Relations of \am-algebras in the tree language.\label{picture}}
\end{figure}
Elements of $\ucalam(n)$ are equivalence classes of linear
combinations of trees from $\treem(n)$, modulo the relations of
Figure~\ref{picture}. Surprisingly, there exists, for any $n\geq 1$, a
very natural basis for the $\bk$-vector space $\ucalam(n)$.

\begin{definition}
\label{90}
We say that a tree $T\in \treem(n)$ is {\em admissible\/} if all its unary
vertices (i.e.~vertices with just one input edge) are input vertices
of the tree $T$.
\end{definition}

Let us denote by $\admm(n) \subset \treem(n)$ the subset consisting of
admissible trees. Then the following proposition holds.

\begin{proposition}
\label{bez_odezvy}
Admissible $n$-trees form, for $n \geq 1$ a linear basis of $\ucalam(n)$,
\[
\ucalam(n) \cong \Span(\admm(n)).
\]
\end{proposition}

\noindent
The proposition follows from the following lemma.
\begin{lemma}
\label{7}
There exists a unique map
$
\Omega : \treem (n) \to \Span(\admm(n))
$
such that, for each $T\in \treem(n)$,
\[
\Omega(T) \equiv T \mbox{ modulo the relations of Figure~\ref{picture}}. 
\]
\end{lemma}

\begin{proof}
If $T$ is admissible, we put
$\Omega(T):= T$. If $T$ is not admissible, it contains a subtree of
the form $\xi_1 \circ_1 \xi_j$, with some $1\leq j \leq m$, or, pictorially,
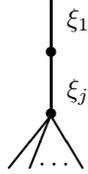
\begin{center}
\setlength{\unitlength}{0.0055in}%
\begin{picture}(70,160)(120,575)
\thicklines
\put(160,735){\line( 0,-1){110}}
\put(160,625){\line(-4,-5){ 40}}
\put(160,625){\line(-2,-5){ 20}}
\put(160,625){\line( 3,-5){ 30}}
\put(160,680){\makebox(0,0)[b]{$\bullet$}}
\put(160,621){\makebox(0,0)[b]{$\bullet$}}
\put(163,573){\makebox(0,0)[b]{$\cdots$}}
\put(175,705){\makebox(0,0)[lb]{{$\xi_1$}}}
\put(175,635){\makebox(0,0)[lb]{{$\xi_j$}}}
\end{picture}
\end{center}
We may also assume that $j \geq 2$, otherwise $T \equiv 0$ modulo $\ucalim$.
Using the relations of Figure~\ref{picture}, we can replace the above
subtree 
by a linear combination of trees of the form
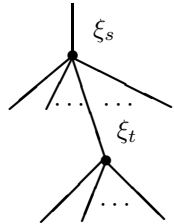
\begin{center}
\unitlength 1.4mm
\thicklines
\begin{picture}(44.50,20.00)(15,15)
\put(35.17,35.00){\line(0,-1){5.00}}
\put(35.17,30.00){\line(-6,-5){6.00}}
\put(35.17,30.00){\line(-1,-2){2.50}}
\put(35.17,29.67){\line(1,-3){3.11}}
\put(38.28,20.33){\line(-1,-1){6.11}}
\put(38.33,20.00){\line(-2,-5){2.27}}
\put(38.33,20.00){\line(1,-1){5.50}}
\put(35.17,29.83){\line(2,-1){9.33}}
\put(39.17,15.83){\makebox(0,0)[cc]{$\cdots$}}
\put(39.50,25.33){\makebox(0,0)[cc]{$\cdots$}}
\put(34.85,25.33){\makebox(0,0)[cc]{$\cdots$}}
\put(35.17,30.00){\makebox(0,0)[cc]{$\bullet$}}
\put(38.33,20.00){\makebox(0,0)[cc]{$\bullet$}}
\put(38.17,32.50){\makebox(0,0)[cc]{$\xi_{s}$}}
\put(40.33,22.83){\makebox(0,0)[cc]{$\xi_t$}}
\end{picture}
\end{center}
with $s\geq 2$. This enables us to move, 
by local replacements, the vertex $\xi_1$
towards the inputs of the tree. Repeating this process sufficiently
many times, we get the requisite presentation $\Omega(T)$. It is
easy to see that this presentation is unique, compare also similar
arguments in the proof of~\cite[Theorem~2.3]{markl:dl}.
\end{proof}

Proposition~\ref{bez_odezvy} immediately implies that
\begin{equation}
\label{eq:2}
\dim_{\bk}(\ucalam(n)) = 2^n \cdot b^m_n,\ \mbox { for } m \leq 2,\ n
\leq 1. 
\end{equation}
where $b^m_n$ is the number of all planar rooted reduced (i.e.~without
unary vertices) $n$-trees, all of whose vertices have arities $ \leq m$.
Clearly, $b^m_n$ equals the number of cells of the $m$-skeleton of
Stasheff's associahedron $K_n$, compare also Example~\ref{0}. Another
kind of information about the size of the non-$\Sigma$-operad $\ucalam$
is given in:

\begin{proposition}
\label{Pejsek}
Let us consider the generating function 
\[
\pphi m(t):= \sum_{n\geq 1} \dim_{\bk}(\ucalam(n)) t^n. 
\]
Then $\pphi 1(t)= 2t$ while, for $m\geq 2$, $\phi(t)=\pphi m(t)$ solves
the equation
\begin{equation}
\label{zajicek_Cmuchalek}
\phi(t)-2t = \phi^2(t)(1+\phi(t)+\cdots+\phi^{m-2}(t)).
\end{equation}
\end{proposition}

\begin{proof}
  Formula~(\ref{zajicek_Cmuchalek}) is probably known, though we were not
  able to find a reference in the literature. We give a proof that can
  be easily modified to obtain the analogous
  formula~(\ref{zajicek_Usacek}) for \lm-algebras.

Let $d_n^m := \mbox{card}(\admm(n))$. By Proposition~\ref{bez_odezvy},
$d_n^m := \dim_{\bk}(\ucalam(n))$. Let us denote, just in this proof,
by $\ttt_n$ the set of all planar reduced $n$-trees all of whose
vertices have at most $m$ input edges.  Let $r_n :=
\mbox{card}(\ttt_n)$. We claim that, for any $n\geq 2$,
\begin{equation}
\label{kohoutek}
r_n =
\sum_{i_1+i_2=n}\hskip -.4em r_{i_1}r_{i_2}+ \hskip -.4em
\sum_{i_1+i_2+i_3=n}\hskip -.6em r_{i_1}r_{i_2}r_{i_3}+
\cdots + \hskip -.7em
\sum_{i_1+\cdots +i_m=n}\hskip -.7em r_{i_1}\cdots r_{i_m}
\end{equation}
and
\begin{equation}
\label{slepicka}
d_n^m  = 2^n r_n.
\end{equation}
The first equation follows from the decomposition
\[
\ttt_n = \ttt_{n,2}\sqcup \ttt_{n,3}\sqcup \cdots \sqcup \ttt_{n,m},
\]
where $\ttt_{n,i}\subset \ttt_n$ is the subset of trees whose root
has $i$ input edges. The second equation follows from the
observation that each admissible tree $T\in \admm(n)$ can be obtained
from a tree $S \in \ttt_n$ by grafting terminal unary trees at the
input legs of $S$; there are $2^n$ ways to do this. If we denote
$\psi(t) := \sum_{n\geq 1}r_nt^n$, then~(\ref{kohoutek}) implies
\[
\psi(t)-t = \psi^2(t)(1+\cdots + \psi^{m-2}(t)),
\]
while~(\ref{slepicka}) implies that $\phi(t)= \psi(2t)$. This finishes
the proof.
\end{proof}

\begin{example}
{\rm\
For $m=\infty$, equation~(\ref{zajicek_Cmuchalek}) can be rewritten as
\[
2\phi^2(t)- \phi(t)(1+2t) + 2t =0.
\]
For $m=2$,~(\ref{zajicek_Cmuchalek}) 
gives $\phi(t)-2t = \phi^2(t)$,
therefore
\[
\phi(t) = \frac{1 - \sqrt{(1-8t)}}2.
\]
A little exercise in Taylor expansions results in
\[
d^2_n = \frac{2^{2n-1}(2n-3)!!}{n!} \hskip 1em \mbox { for } n \geq 2,
\]
where as usual 
\[
(2n-3)!! = 1 \cdot 3 \cdot 5 \cdots (2n-3).
\]
}\end{example}

\section{Free \am-algebras}

The following theorem relies on the notation introduced in the previous
section. Recall namely that 
$\admm(n)$ denotes the set of admissible $n$-trees in
the sense of Definition~\ref{90}, $\ucalam(n)\cong \Span(\admm(n))$ and
$\Omega$ is the map of Lemma~\ref{7}.
Recall also that we denoted by $\xi_j$ the $j$-corolla.
If $T_1,\ldots,T_j$ are trees, then
$\xi_j(T_1,\ldots,T_j)$ is the tree obtained by
grafting the tree 
$T_i$ at the $i$-th input edge of the tree $\xi_j$, $1\leq
i \leq j$.

\begin{theorem}
Let $X$ be a graded vector space. Then the free \am-algebra  $\ \frm(X)$
on $X$ can be described as \ $\frm(X) = \bigoplus_{n\geq 1}\frm{}^n(X)$, with
\[
\frm{}^n(X) := \ucalam(n)\otimes X^{\otimes n}.
\]
The structure operations $\mu_1,\ldots,\mu_m$ are defined as follows. Let
$\VECTOR {v}_j \in X^{\otimes k_j}$ and $T_j \in \admm(k_j)$, $1\leq
j\leq k$. Then 
\[
\mu_k(T_1\otimes \VECTOR v_1,\ldots,T_k\otimes \VECTOR v_k) :=
\Omega(\xi_k(T_1,\ldots,T_k)) \otimes (\VECTOR v_1\otimes \cdots
\otimes \VECTOR v_k).
\]
\end{theorem}

\begin{proof}
  It follows from general
  theory~\cite[Section~II.1.4]{markl-shnider-stasheff:book} that the free
  $\underline{\EuScript P}$-algebra $F_{\underline{\EuScript P}}(X)$
  on a graded vector space $X$, where $\underline{\EuScript P}$ is a
  non-$\Sigma$ operad, is given by $F_{\underline{\EuScript P}}(X) =
  \bigoplus_{n\geq 1}F_{\underline{\EuScript P}}^n(X)$, with
  $F_{\underline{\EuScript P}}^n(X):= \underline{\EuScript
    P}(n)\otimes X^{\otimes n}$, with the algebra structure induced
  from the operad structure of $\underline{\EuScript P}$. Our theorem
  is then obtained by taking $\underline{\EuScript P} := \ucalam$ and
  applying the results on the structure of this operad proved in the
  previous section.
\end{proof}

Observe that the free algebra $\frm (X)$ is bigraded. The {\em internal\/}
grading is the grading of the underlying vector space, while
the {\em external\/} grading is given by the decomposition
\[
\frm(X) = \bigoplus_{n\geq 1}\frm{}^n(X).
\] 
Because all main structures considered in this note (vector
spaces, \am-algebras, \lm-algebras, etc.) have an implicit internal
grading, we will usually mean by saying that an object is {\em graded\/} the 
presence of an external grading. We believe that the actual meaning
will always be clear from the context.

\begin{example}
\label{0}
{\rm
  There exists a simple way to encode elements of
  $\admm(n)$.  Let $\bm(n)$ be the set of all meaningful bracketings
  of strings $\epsilon_1\epsilon_2\cdots\epsilon_n$, where
  $\epsilon_i\in \{\circ,\bullet\}$ and no pair of brackets
  encompasses more than $m$ terms. A moment's reflection shows that
  $\admm(n) \cong \bm(n)$, $n\geq 1$.  Thus
\begin{eqnarray*}
\adm m(1) &=& \{\c,\b\},\ m\geq 1,
\\
\adm 1(2) &=& \emptyset, \ 
\admm(2) = \{(\c\c), (\c\b), (\b\c), (\b\b) \},\ m\geq 2,
\\
\adm1(3)& =& \emptyset,
\\
\adm2(3)&=&\{((\c\c)\c),((\b\c)\c),((\b\b)\c),
\ldots,(\c(\c\c)),(\b(\c\c)),(\b(\b\c)),\ldots\},
\\
\adm m(3)&=& 
\hskip .2em \adm2(3) \cup \{(\c\c\c),(\b\c\c),(\b\b\c),\ldots \},\  m\geq 3,
\mbox{ etc.}
\end{eqnarray*}

We conclude that $\frm(X)$ is spanned by all meaningful bracketings of
strings of ``variables'' $x$ and $\pa x$, $\deg(\pa x) = \deg(x)-1$,
$x \in X$, with no pair of brackets encompassing more than $m$ terms.
Thus $\frmfff 1(X) = \frmfff 1{}^1(X)$ is the free differential
complex generated by $X$, $\frmfff 1(X) \cong X\oplus\pa X$. Similarly,
\begin{eqnarray*}
\frm{}^2 (X) &\cong & 
(X \otimes X)\oplus (\pa X \otimes X)\oplus 
(X \otimes \pa X)\oplus (\pa X \otimes \pa X),\ \mbox{for $m \geq 2$,}
\\
\frmfff 2{}^3 (X) & \cong& 
((X \otimes X ) \otimes X) \oplus 
((\partial X \otimes X ) \otimes X) \oplus 
((X \otimes \pa X) \otimes X)
\\
&&\hskip 1.5em \oplus\ ((X \otimes X) \otimes \pa X)
\oplus  
((\pa X \otimes \pa X) \otimes X) 
\oplus\ ((\pa X \otimes X) \otimes\pa X)
\\
&&\hskip 1.5em  \oplus\ ((X \otimes \pa X) \otimes \pa X)
 \oplus ((\pa X \otimes \pa X) \otimes \pa X)
\oplus (X \otimes (X  \otimes X))
\\
&&\hskip 1.5em  \oplus\ 
(\partial X \otimes( X \otimes X)) \oplus 
(X \otimes (\pa X \otimes X)) \oplus (X \otimes( X \otimes \pa X))
\\
&&\hskip 1.5em 
\oplus\  (\pa X \otimes( \pa X \otimes X)) 
\oplus(\pa X \otimes 
(X \otimes\pa X)) \oplus (X \otimes( \pa X \otimes \pa X))
\\
&&\hskip 1.5em  \oplus\ (\pa X \otimes (\pa X \otimes \pa X)),
\\
{\frm}{}^3(X)  &\cong& 
\frmfff 2{}^3 (X) \oplus (X \ot X \ot X)
 \oplus (\pa X \ot X \ot X)  \oplus (X \ot\pa  X \ot X)  
\oplus (X \ot X \ot
 \pa X)
\\
&&\hskip 1.5em  \oplus\ (\pa X \ot\pa  X \ot X) \oplus (\pa X \ot X \ot\pa  X) 
\oplus (X \ot \pa X \ot
\pa X) 
\\
&&\hskip 1.5em 
\oplus\ (\pa X \ot\pa  X \ot\pa  X),\ \mbox{for $m \geq 3$, etc.}
\end{eqnarray*} 
We see that
\[
\dim(\hskip .2em \frm{}^n(X)) = (2 \dim(X))^n \cdot b^m_n,\ n \geq 1,\
m \geq 2,
\]
where $b^m_n$ is as in~(\ref{eq:2}).
}
\end{example}

\noindent
{\bf Unital {\em {\rm \am}}-algebras.}
Recall~\cite[page~148]{markl:JPAA92} that an \am-algebra
  $A=(A,\{\mu_n\})$ is {\em unital\/} if there exists an element $1=
  1_A \in A$ such that $\mu_2(a,1_A)= \mu_2(1_A,a) = a$, for all $a\in
  A$, and that $\mu_n(\rada a1{i-1},1_A,\rada a{i+1}n)= 0$, for all $n
  \not= 2$ and $1\leq i \leq n$. The {\em free unital \am-algebra\/}
  can be constructed as $\Span(1)\oplus \frm (X)$, with the structure
  operations defined in an obvious way.

\vskip 2mm
\noindent 
{\bf Important modifications.}  One may as well consider, for any $n <
m$, the forgetful functor $\Box : \mbox{\catam} \to \mbox {\catan}$.
This functor has a left adjoint $\frmff mn : \mbox{\catan} \to \mbox
{\catam}$ which can be described as follows. Given an \an-algebra
$(V,m_1,m_2,\ldots,m_n)$, consider the free \am-algebra $\frm(V)$ on
the graded vector space $V$, and let $\mu_1,\ldots,\mu_m$ be the
structure operations of $\frm(V)$.  Then
\[
\frmff mn(V,m_1,m_2,\ldots,m_n) \cong \frm(V)/(\mu_i = m_i,\ i \leq n).
\] 
A very important special case is $n=1$ when we obtain a left adjoint
$\frmff m1 : \mbox {\catdgVect} \to \mbox {\catam}$ to the forgetful
functor $\Box :\mbox {\catam} \to \mbox {\catdgVect}$ to the category
of {\em differential\/} graded vector spaces. One can easily verify
that
\begin{equation}
\label{11}
\frmff m1 (V,\partial) \cong \bigoplus_{n \geq 1} (\ucalamdg(n),\partial)
\otimes (V,\partial)^{\otimes n},
\end{equation}
where $(\ucalamdg,\partial)$ is now the non-$\Sigma$ operad describing
\am-algebras as algebras in the monoidal category of {\em
  differential\/} vector spaces. Let us
recall~\cite[page~1493]{markl:zebrulka} that $\ucalamdg =
\uGamma(\xi_2,\ldots, \xi_m)$ with the differential
\[
d(\xi_n) := \sum (-1)^{(b+1)(i+1)+b}\cdot  \xi_a \circ_i \xi_b,
\]
where the summation runs over all $a,b \geq 2$ with $a+b = n+1$ and $1
\leq i \leq a$.
Formula~(\ref{11}) is in a perfect harmony
with~\cite[Definition~II.1.24]{markl-shnider-stasheff:book} describing
free operad algebras in a general monoidal category.

\begin{example}
{\rm
Clearly $\frmff m1{}^1(V,\pa) \cong (V,\pa)$, and,
for $m \geq 2$, $\frmff m1{}^2 (V,\pa) \cong (V \otimes V)$ with the
induced differential. As graded vector spaces 
\begin{eqnarray*}
\frmff 21{}^3 (V,\pa)  & \cong & 
((V \otimes V ) \otimes V) \oplus  (V \otimes (V  \otimes V)),
\\
\frmff m1{}^3(V,\pa) & \cong &\frmff 21{}^3 (V) 
\oplus (V \ot V \ot V), \mbox { for $m \geq 3$,}
\\
\frmff 21{}^4(V,\pa)  
& \cong & (((V \ot V) \ot V) \ot V) \oplus
(V \ot (V \ot (V \ot V))) \oplus ((V \ot V) \ot (V \ot V))
\\
&& \hskip 2em
\oplus\ ((V \ot (V \ot V)) \ot V) \oplus (V \ot ((V \ot V) \ot V)),
\\
\frmff 31{}^4(V,\pa) & \cong & 
\frmff 21{}^4(V,\pa) \oplus
((V \ot V \ot V) \ot V) \oplus (V \ot (V \ot V \ot V)) 
\oplus ((V \ot V) \ot V \ot V) 
\\
&&  \hskip 2em
\oplus\  (V \ot (V \ot V) \ot V) \oplus (V \ot V \ot (V \ot V)),
\\
\frmff m1{}^4(V,\pa) & \cong & 
 \frmff 31{}^4(V,\pa) \oplus 
(V \ot V \ot V \ot V \ot V) \mbox { for $m \geq 4$, etc.} 
\end{eqnarray*} 
The action of the differential in $\frmff m1(V,\pa)$ can be 
easily read off from the axioms of \am-algebras, for example
\begin{eqnarray*}
\lefteqn{
\pa(x \ot y \ot z) =  ((x \ot y) \ot z)- (x \ot (y \ot z))\hskip 10em }
\\
&&\hskip 10em
- (\pa x \ot y \ot z) -(-1)^{|a|} (x \ot \pa y \ot z) 
-(-1)^{|a| + |b|} (x \ot y \ot \pa z),
\end{eqnarray*} 
where $(x \ot y \ot z) \in \frmff m1{}^3(V)$, $m \geq 3$.  The above
formula illustrates a surprising fact that the differential in $\frmff
m1 (V,\pa)$ is nontrivial even if the differential $\pa$ on $V$ is
zero, because in this case still
\[
\pa(x \ot y \ot z) =  ((x \ot y) \ot z)- (x \ot (y \ot z))\ !
\]
The dimension of the homogeneous components of $\frmff m1(V,\pa
V)$ clearly equals
\[
\dim{\frmff m1{}^n(V,\pa)} = (\dim(V))^n \cdot b^m_n, \mbox{ for }
n \geq 1,\ m \geq 2,
\]
where the integers $b^m_n$ are as in~(\ref{eq:2}).
}\end{example}

\section{Free strongly homotopy Lie algebras}
\label{to_jsem_zvedav_jestli_aspon_neco_bude_prijato}

In this section we indicate how to modify previous results to strongly
homotopy Lie algebras. Recall (\cite{lada-stasheff:IJTP93}
or~\cite[Definition~2.1]{lada-markl:CommAlg95}) that an \LLL
m-structure on a graded vector space $L$ is a system $\{l_k;\ 1\leq
k\leq m,\ k<\infty\}$ of degree $k-2$ linear maps $l_k:\otimes^k L \to
L$ which are antisymmetric in the sense that
\begin{equation}
\label{antisymmetry}
l_k(\rada x{\sigma(1)}{\sigma(k)})=\chi(\sigma)l_k(\rada x1n)
\end{equation}
for all $\sigma \in \Sigma_n$ and $\rada x1n\in L$.
Here $\chi(\sigma):= \sgn(\sigma)\cdot \epsilon(\sigma)$, where
$\epsilon(\sigma)$ is the Koszul sign of the permutation $\sigma$, 
see~\cite[page~2148]{lada-markl:CommAlg95} for details.
Moreover, the
following generalized form of the Jacobi identity is supposed to be
satisfied for any $n\leq m$:
\begin{equation}
\label{Jacobi}
\sum_{i+j=n+1}\sum_\sigma
\chi(\sigma)(-1)^{i(j-1)}l_j(l_i(\rada x{\sigma(1)}{\sigma(i)}),\rada
x{\sigma(i+1)}{\sigma (n)})=0,
\end{equation}
where the summation is taken over all $(i,n-i)$-unshuffles with
$i\geq 1$.

The operad $\callm = \{\callm(n)\}_{n \geq 1}$ describing \LLL
m-algebras as algebras in the category of graded vector spaces can be
constructed as follows.  Let $\calfm = \{\calfm(n)\}_{n\geq 1}$ be the
free graded $\Sigma$-operad $\Gamma(\zeta_1,\zeta_2,\ldots,\zeta_m)$
generated by the operations $\zeta_i$, $l\leq i \leq m$, where each
$\zeta_i$ has $i$ inputs, degree $i-2$ and spans the signum
representation of the symmetric group $\Sigma_i$.  Using trees, this
operad can be described as follows. Let $\treem(n)$ be now the set of
all directed labeled (non-planar, abstract) trees with $n$-input
edges, all of whose vertices have at most $m$ input edges. We again
admit also `unary' vertices with only one input edge and `labeled'
means that the input legs are labeled by natural numbers $1,\ldots,n$.
Then, as before,
\[
\calfm(n) \cong \Span(\treem(n)).
\]
We put $\callm := \calfm/\calim$, where the operadic ideal $\calim$ is now
generated by relations whose pictorial presentation is indicated
in Figure~\ref{liquid}.

\begin{figure}
\begin{center}
\setlength{\unitlength}{0.007in}
\begin{picture}(314,195)(-95,560)
\thicklines
\put(60,645){\makebox(0,0)[cc]{$0=\displaystyle\sum_{i+j=n+1}\sum_\sigma
\chi(\sigma)(-1)^{i(j-1)}\cdot%
\left(\rule{0mm}{20mm}\hskip55mm \right)$}}
\put(14,0){
\put(160,730){\line( 0,-1){ 60}}
\put(160,670){\line(-5,-6){ 50}}
\put(110,610){\line(-1,-1){ 50}}
\put(110,610){\line(-2,-5){ 20}}
\put(110,610){\line( 4,-5){ 40}}
\put(160,670){\line( 1,-4){ 25}}
\put(185,570){\line( 0,-1){  5}}
\put(160,670){\line( 3,-2){150}}
\put(145,660){\makebox(0,0)[rb]{$\xi_j$}}
\put(98,610){\makebox(0,0)[rb]{$\xi_i$}}
\put(160,665){\makebox(0,0)[b]{$\bullet$}}
\put(110,607){\makebox(0,0)[b]{$\bullet$}}
\put(115,570){\makebox(0,0)[b]{$\cdots$}}
\put(0,-6){
\put(235,580){\makebox(0,0)[b]{$\cdots$}}
\put(53,540){\makebox(0,0)[b]{$\sigma(1)$}}
\put(92,540){\makebox(0,0)[b]{$\sigma(2)$}}
\put(155,540){\makebox(0,0)[b]{$\sigma(i)$}}
\put(207,550){\makebox(0,0)[b]{$\sigma(i+1)$}}
\put(320,550){\makebox(0,0)[b]{$\sigma(n)$}}
}
}
\end{picture}
\end{center}
\caption{Relations for strongly homotopy Lie algebras in the
tree language.\label{liquid}}
\end{figure}

By similar arguments as in the previous sections we easily infer that
there is a natural basis of the $\bk$-vector space $\callm(n)$, formed
by the set $\admm(n)$ of admissible $n$-trees, whose definition is
formally the same as in Definition~\ref{90}, only removing everywhere
the adjective ``planar.'' The size of the operad $\callm$ is
described in the following proposition whose proof is similar to that
of Proposition~\ref{Pejsek}.

\begin{proposition}
\label{jsem_nachlazeny}
The generating function  
$\varomega(t) =\vvaromega m(t):= 
\sum_{n\geq 1} \frac1{n!}\dim_{\bk}(\callm(n)) t^n$
satisfies
the equation
\begin{equation}
\label{zajicek_Usacek}
\varomega(t)-2t = \varomega^2(t)\left(\frac1{2!}+\frac1{3!}\varomega(t)+
\cdots+\frac1{m!}\varomega^{m-2}(t)\right).
\end{equation}
\end{proposition}

\begin{example}
  {\rm\ For $m=\infty$, equation~(\ref{zajicek_Usacek}) becomes
    $2\varomega(t)-2t-1 = \exp(\varomega(t))$.  For $m=2$,
    equation~(\ref{zajicek_Usacek}) says that $\varomega(t)=
    1-\sqrt{1-4t^2}$, therefore 
\[
\dim_{\bk}(\calldve(n)) = 2^n \cdot (2n-3)!!,\ n\geq 1.
\]
}\end{example}

As in Lemma~\ref{7}, we have a map $\Omega : \treem (n) \to
\Span(\admm(n))$ with the property that $\Omega(T) \equiv T$ modulo
relations of Figure~\ref{liquid}. Free \LLL m-algebras are then
described as follows.

\begin{theorem}
Let $Y$ be a graded vector space. The free \LLL m-algebra $\frml(Y)$
generated by $Y$ can defined as $\frml(Y) = 
\bigoplus_{n\geq 1}\frml{}^n(Y)$, where
\[
\frml{}^n(Y) := \callm(n)\otimes_{\Sigma_n} Y^{\otimes n}.
\]
The structure operations $l_1,\ldots,l_m$ are given by the 
following rule. 
Let $\VECTOR {v}_j \in Y^{\otimes k_j}$ and $T_j \in \admm(k_j)$, $1\leq
j\leq k$. Then 
\[
l_k(T_1\otimes_{\Sigma_{k_1}} 
\VECTOR v_1,\ldots,T_k\otimes_{\Sigma_{k_1}} \VECTOR v_k) :=
\Omega(\xi_k(T_1,\ldots,T_k)) \otimes_{\Sigma_{k_1+\cdots+k_k}} 
(\VECTOR v_1\otimes \cdots
\otimes \VECTOR v_k).
\]
\end{theorem}

\begin{remark}
  {\rm In exactly the same fashion as for \am-algebras one may define
    left adjoints $\callff mn : \mbox{\catln} \to \mbox{\catlm}$ to
    the forgetful functor $\Box : \mbox{\catlm} \to \mbox{\catln}$, $n
    < m$, as
\[
\callff mn\ (L,\lambda_1,\lambda_2,\ldots,\lambda_n)
=
\frml (L)/(\lambda_i = l_i,\ i \leq n),
\] 
where $l_i$ are the structure operations of $\frml(L)$.
}
\end{remark}

\section{Reflections on the PBW theorem}

As in the case of ordinary Lie algebras, there exists the
symmetrization functor from the category \catam\ of \am-algebras to the
category \catlm\ of \lm-algebras. More precisely, we proved 
in~\cite[Theorem~3.1]{lada-markl:CommAlg95}:

\begin{theorem}
\label{Pejsek_s_Kocickou}
Any \am-structure $\lbrace \mu_n:\bigotimes ^nV\longrightarrow
V\rbrace$
on a differential graded vector space $V$ induces an \lm-structure
$\lbrace l_n:\bigotimes ^nV \longrightarrow V \rbrace$ on the same
differential graded vector space, with
\[
l_n(v_1\otimes \dots \otimes v_n):=\sum_{\sigma\in \cals_n}\chi
(\sigma )
\mu_n(v_{\sigma (1)}\otimes \dots \otimes v_{\sigma (n)}),\ 2\leq
n\leq m.
\]
This correspondence defines a functor (the {\em symmetrization\/}) 
$(-)_L:{\tt A}(m) \longrightarrow {\tt L}(m)$.

\end{theorem}

In~\cite[Theorem~3.3]{lada-markl:CommAlg95}, we also proved that: 

\begin{theorem}
  There exists a functor $\mathcal U_m : {\tt L}(m) \longrightarrow
  {\tt A}(m)$ that is left adjoint to the functor $(-)_L:{\tt A}(m)
  \longrightarrow {\tt L}(m)$.
\end{theorem}

The algebra $\mathcal U_m(L)$ of the previous theorem is called
the {\em universal enveloping\/} \am-{\em algebra} 
of the \lm-algebra $L$.
We gave, in~\cite[page~2154]{lada-markl:CommAlg95},
the following explicit construction of
$\calu_m(L)$. 

Start with the free {\em unital\/} $\AM$-algebra $\frm(L)$ generated by the
underlying vector space $L$ of the \lm-algebra $L=(L,\{l_n\})$ and let
$\lbrace \mu _n\rbrace$ be the $\AM$-structure maps of $\frm(L)$. Let
$I$ denote the \am-ideal in $\frm(L)$ generated by the relations
\[
\sum_{\sigma \in S_n} \chi (\sigma) \mu_n (\xi _{\sigma (1)},\dots ,
\xi _{\sigma(n)})
= l_n(\xi _1, \dots , \xi_n),\mbox{ for } \xi_1,\dots,\xi_n\in L,\ n
\leq m.
\]
Then we put $\mathcal U_m(L) := \frm(L) / I$. This universal
enveloping \am-algebra is equipped with the canonical map (in fact, an
inclusion) $\iota: L \to \mathcal U_m (L)$.

\begin{remark}
  {\rm Hinich and Schechtman introduced
    in~\cite{hinich-schechtman:AdvinSovMath1993}, for any \lm-algebra
    $L$ (they did it for $m=\infty$, but the general case is an
    obvious modification), an {\em associative\/} algebra (which they
    also called the universal enveloping algebra), having the property
    that the category of modules (in a suitable sense, see for
    example~\cite[Definition~5.1]{lada-markl:CommAlg95}) over an
    \lm-algebra $L$ is equivalent to the category of left modules over
    the universal enveloping associative algebra of $L$.
  
    It is unclear whether there exists a functor ${\EuScript G} :
    \mbox {\catlm} \to \mbox {\catam}$ such that the category of
    modules over an \lm-algebra $L$ would be naturally equivalent to
    the category of left modules, in the sense
    of~\cite[page~157]{markl:JPAA92}, over the \am-algebra $\EuScript
    G (L)$. Our universal enveloping \am-algebra functor whose
    definition we recalled above probably does not have this property.
    While it can be easily shown that each left module over the
    \am-algebra\ $\mathcal U_m(L)$ naturally defines a module over
    $L$, we doubt that each $L$-module uniquely extends into a left
    module over $\mathcal U_m(L)$.  The classical universal enveloping
    associative algebra of a strict Lie algebra is constructed as a
    left adjoint to the symmetrization functor. By a miracle, this
    functor also happens to describe the category of $L$-modules, but this
    mystic phenomenon does not take place in the category of
    \lm-algebras.
  
    Moreover, there is no reason to believe that such a functor
    $\EuScript G$ exits. The category of modules over a given
    \lm-algebra $L$ is a well-behaved abelian category, therefore it
    should be the category of modules over an {\em associative\/}
    algebra, which turns out to be the one constructed by Hinich and
    Schechtman. On the other hand, we do not see why this abelian
    category should be at the same time the category of left modules
    over some \am-algebra.  Let us close this remark with

\begin{problem}
  Let $L$ be an {\em \lm\/}-algebra.  Does there exist an {\em
    \am\/}-algebra $A$ with the property that the category of modules
  over $L$ is equivalent to the category of left modules over $A$?
\end{problem}
}\end{remark}

The simplest form of the P(oincar\'e)-B(irkhoff)-W(itt) theorem for an
ordinary Lie algebra $L= (L,[-,-])$ says that the associated graded
algebra $G^*(L)$ of the universal enveloping algebra ${\mathcal U}(L)$
of $L$ is isomorphic to the free commutative associative algebra
$S^*(L)$ (the polynomial ring) on the vector space $L$.  More
precisely, recall that
\[
\calu(L)= T(L)/(x\ot y - y \ot x = [x,y],\ x,y \in L),
\] 
where $T(L)$ is the free associative algebra (tensor algebra) on the
vector space $L$.  The algebra $\ul$ is filtered, the ascending
filtration being given by $\calu_p(L) :=$ vector space generated by
linear combinations of elements which can be presented as a product of
$\ \leq p$ elements in the augmentation ideal of $T(L)$.  Clearly
$\calu_{0}(L)= \bk$ and $\calu_1(L) = \bk \oplus L$. If $G^*(L)$ denotes
the associated graded algebra, then one proves that the natural map
$L=G^1(L) \to G^*(L)$ is a monomorphism and $G^*(L)$ is generated by
$L$.  It is immediate to see that $xy = yx$, for arbitrary $x,y \in L
\subset G^*(L)$. We formulate the following trivial observation.

\begin{observation}
\label{opicak_Fuk}
Suppose that an associative algebra $A$ admits a set of mutually
commutative  generators. Then $A$ is commutative. 
\end{observation}

{}From this observation we infer that $G^*(L)$ is commutative,
therefore there exists, by the universal property, a commutative
algebra map $S^*(L)\to G^*(L)$, induced by the inclusion $L \subset
G^*(L)$. The PBW\footnote{not to be mistaken with
  Praha-Berlin-Warszawa cycling race} theorem says that this map is an
isomorphism of graded commutative associative algebras.

Let us try to guess which form the PBW theorem for the universal
enveloping \am-algebra $\calu_m(L)$ of an \lm-algebra $L$ may possibly
have. The classical PBW theorem compares the associated graded of
\hskip .1em $\ul$ to a commutative associative algebra.  Commutative
associative algebras occur in the symbolic exact sequence (which may
be found in~\cite[page~228]{ginzburg-kapranov:DMJ94})
\[
\mbox{Commutative associative algebras} \subset \mbox{Associative  algebras}
\stackrel{\pi}{\longrightarrow} \mbox{Lie algebras.} 
\]
Here the `projection' $\pi$ is given by the symmetrization of the
associative product. Commutative associative 
algebras are also the quadratic duals of Lie algebras 
(and vice versa)~\cite[Theorem~2.1.11]{ginzburg-kapranov:DMJ94}, 
but this coincidence is misleading. Summing up,
the classical PBW theorem says that the associated graded $G^*(L)$
is the free algebra in the kernel of the symmetrization map~$\pi$.

One would naturally expect something similar 
also in the \lm-algebra case. The kernel
of the symmetrization map $(-)_L$ of 
Theorem~\ref{Pejsek_s_Kocickou} clearly
consists of \am-algebras
$A = (A,\{\mu_n\})$ such that
\begin{equation}
\label{Ferdik}
\sum_{\sigma \in S_n} \chi(\sigma)
\mu_n(a_{\sigma(1)},\ldots,a_{\sigma(n)})  = 0,
\ n\geq 2,\ \rada a1n\in A.
\end{equation}
Therefore one tends to believe that \am-algebras enjoying this form of
symmetry are analogs of the polynomial algebra $S^*(L)$ from the
classical PBW theorem. This is, however, not exactly so.

To see why, let $L = (L,\{l_n\})$ be an \lm-algebra and define
inductively an ascending filtration of the universal enveloping
\am-algebra $\calu_m(L)$ by $\calu_{m,0}(L):=\bk$, $\calu_{m,1}(L) :=
\calu_{m,0}\cup \iota(L)$ and $\calu_{m,p}(L)$ being the vector space
generated by elements of the form $\mu_k(\rada x1k)$, $k\geq 2$, $x_i
\in \calu_{m,p_i}(L)$ for $1\leq i\leq k$ and $\sum_1^k p_i \leq p$.
The associated graded \am-algebra $G^*_m(L)$ is defined by
\[
G^*_m(L) = \bigoplus_{q\geq 0}G^q_m(L),\
G^q_m(L) :=\calu_{m,q}(L)/ \calu_{m,q-1}(L)
\]
with the \am-structure maps induced by that of $\calu_m(L)$.
It is easily seen that $G^*_m(L)$ 
is generated by the image of the canonical
inclusion $\iota: L \hookrightarrow \calu_{m,1}(L) \to G^*_m(L)$. 
A moment's reflection shows that the elements of $L\subset G_m^*(L)$
`mutually commute' in the sense that~(\ref{Ferdik}) is satisfied for
$a_1,\ldots,a_n \in L$. Surprisingly, 
this {\em does not\/} imply that~(\ref{Ferdik})
is satisfied for {\em all elements of $G^*_m(L)$\/} -- 
there is no analog of
Observation~\ref{opicak_Fuk} for \am-algebras!
One must instead consider \am-algebras $S_m(V,\pa)$
defined as
\[
S_m(V,\pa):= \frmff m1(V,\pa)/{\mathfrak I}_m,
\]
where ${\mathfrak I}_m$ is the ideal generated by the relations
\[
\sum_{\sigma \in S_n} \chi(\sigma)
\mu_n(v_{\sigma(1)},\ldots,v_{\sigma(n)})  = 0\ \mbox { for } 
\ 2 \leq n \leq m\ \mbox { and }\ \rada v1n\in V.
\]
The algebra $S_m(V,\pa)$ is (externally) graded, with the grading induced from
the (external) grading of the free \am-algebra $\frmff m1(V,\pa)$. 
By a miracle, $S_m^*(V,\pa)$ has the 
following very explicit description.

Choose a basis $\{f_i\}_{i\in I}$ of $V$ indexed by a linearly ordered
set $I$. Let $\treem{}^n(V)$ be the set of all planar rooted reduced
$n$-trees $S$ whose vertices have arity $\leq m$. We moreover
suppose that the input legs of $S$ are labeled by elements of the
basis $\{f_i\}_{i\in I}$ and that the labels satisfy the following
condition:

Let $v$ be an input vertex of $S$ and let $\rada e1k$ be the
input legs of $v$, written in the order induced by the imbedding of $S$
into the plane. Let $\rada f{i_1}{i_k}$ be the corresponding labels.
Then it is \underline{not} true that ${i_1}\geq i_2 \geq \cdots \geq
i_k$.

\begin{proposition}
For each $n\geq 1$, $S^n_m(V) \cong \Span(\treem{}^n(V))$.
\end{proposition}

\begin{proof}
  It follows from~(\ref{11}) that $\frmff m1(V,\pa)$ is spanned by the
  set of planar directed reduced $n$-trees with all vertices of arity
  $\leq m$ and the input legs labeled by the chosen basis
  $\{f_i\}_{i \in I}$ of $V$. Let $v$ be an input vertex of $S$,
  $e_1,\ldots,e_k$ input legs of $v$ written in the order induced by
  the plane and $\rada f{i_1}{i_k}$ their labels.
  
  Suppose ${i_1}\geq i_2 \geq \cdots \geq i_k$.  If $i_1 = i_k$, then
  clearly $S \equiv 0$ modulo ${\mathfrak I}_m$.  If $i_1 \not =
  i_k$, we may replace $S$ modulo ${\mathfrak I}_m$ by a sum of trees
  that differ from $S$ only by labels of $\rada e1k$ which are
  changed to $\rada i{\sigma(1)}{\sigma(k)}$, for some $\sigma \in
  \Sigma_k$ such that it is not true that
  ${i_{\sigma(1)}}\geq i_{\sigma(2)} \geq \cdots \geq i_{\sigma(k)}$.
  Repeating this process at any input vertex of $S$, we replace $S$
  modulo ${\mathfrak I}_m$ by a sum of trees from $\treem{}^n(V)$.
  This replacement is obviously unique.
\end{proof}

Let us formulate our version of the PBW theorem for the universal
enveloping \am-algebra of an \lm-algebra.

\begin{theorem}
The canonical map $\iota :L \hookrightarrow \calu_{m,1}(L)\to G^*_m(L)$
induces an isomorphism
\[
\rho: S^*_m(L) \cong G^*_m(L)
\]
of (externally) graded \am-algebras.
\end{theorem}

\begin{proof}
  The proof relies on the fact that $\{\rho(S);\ S\in \treem{}^n(V)\}$
  forms a basis of $G^n_m(V)$, which easily follows from the above
  analysis. We leave the details to the reader.
\end{proof}
 
\begin{example}{\rm
Let us analyze the case $m = 2$ when
\lm-algebras are just differential graded
anti-commutative non-associative algebras. For an L(2)-algebra $L =
(L,[-,-],\pa)$, $\frmff 21(L)$ is the free
non-associative dg algebra generated by $(L,\pa)$, the universal enveloping
A(2)-algebra of $L$ is given by
\[
\mathcal U_2(L) = \frmff 21 / (xy - (-1)^{|x||y|}yx = [x,y],\  x,y \in L),
\]
and the associated graded $G^*_2(L)$ is the quotient of $\frmff 21(L)$
by the ideal generated by $xy = (-1)^{|x||y|} yx$, for $x,y \in L$.  }
\end{example}

\begin{remark}{\em
    An ordinary strict Lie algebra $L = (L,[-,-])$ can be considered,
    for any $m \geq 2$, as an \lm-algebra with $l_2 = [-,-]$ and $l_n
    = 0$ if $n \not= 2$. Its universal enveloping \am-algebra
    ${\mathcal U}_m(L)$ is a fully-fledged \am-algebra, therefore an
    object completely different from the ordinary universal enveloping
    associative algebra $\ul$.  Thus our $\mathcal U_m(L)$ cannot be
    considered as a generalization of the classical universal
    enveloping algebra functor.  }
\end{remark}

\section*{Appendix: the tree language}

The following definitions can be found for example
in~\cite[Section~II.1.5]{markl-shnider-stasheff:book}.  By a {\em
  tree\/} we mean a connected graph $T$ without loops. Let $\vert(T)$
denote the set of vertices of $T$. A {\em valence\/} of a vertex $v
\in \vert(T)$ is the number of edges adjacent to $v$. A {\em leg\/} of
$T$ is an edge adjacent to a vertex of valence one, other edges of $T$
are {\em interior\/}.  We in fact discard vertices of valence one at
the endpoints of the legs, therefore the legs are ``half-edges'' having
only one vertex. This can be formalized by introducing (generalized)
graphs as certain sets with involutions,
see~\cite[Definition~II.5.37]{markl-shnider-stasheff:book}.

By a {\em rooted\/} or {\em directed\/} tree we mean a tree with a
distinguished {\em output\/} leg called the {\em root\/}.  The
remaining legs are called the {\em input legs\/} of the tree. A tree
with $n$ input legs is called an {\em $n$-tree\/}.  A rooted tree is
automatically {\em oriented\/}, each edge pointing towards the root.
The edges pointing towards a given vertex $v$ are called the {\em
  input edges\/} of $v$, the number of these input edges is then the
{\em arity\/} of $v$. Vertices of arity one are called unary, vertices
of arity two binary, vertices of arity three ternary,~etc.  A tree is
{\em reduced\/} if it has no unary vertices.  A vertex $v \in
\vert(T)$ is an {\em input vertex\/} of $T$ if all its input edges are
input legs of $T$.  With the above conventions accepted, we must admit
also the unique reduced 1-tree which consist of one edge which is both
its root and the input leg, and no vertices (soul without body).

Each tree can be imbedded into the plane. By a {\em planar tree\/} we
mean a tree with a specified (isotopy class of) imbedding.  Reading
off counterclockwise the input legs of a planar tree, starting from
the leftmost one, gives a linear {\em order} of input legs. In the same
manner, input edges of each vertex of a planar tree are linearly
oriented. Vice versa, specifying linear orders of input edges of each
vertex of a rooted tree $T$ uniquely determines an isotopy class of
imbedding of $T$ into the plane.

\begin{example}{\rm
\begin{figure}[t]
  \centering
\unitlength .5em
\thicklines
\begin{picture}(38,27)(0,2)
\put(19,3){\vector(0,1){25}}
\put(19,22){\line(-1,-1){19}}
\put(19,22){\line(1,-1){19}}
\put(35,6){\line(-1,-1){3}}
\put(3,6){\line(1,-1){3}}
\put(3,6){\line(0,-1){3}}
\put(19,16){\line(-2,-3){6}}
\put(13,7){\line(0,1){0}}
\put(13,7){\line(-1,-1){4}}
\put(9,3){\line(0,1){0}}
\put(13,7){\line(0,-1){4}}
\put(13,7){\line(1,-1){4}}
\put(19,16){\line(5,-6){5}}
\put(24,10){\line(0,1){0}}
\put(24,10){\line(0,-1){7}}
\put(19,22){\makebox(0,0)[cc]{$\bullet$}}
\put(19,16){\makebox(0,0)[cc]{$\bullet$}}
\put(3,6){\makebox(0,0)[cc]{$\bullet$}}
\put(13,7){\makebox(0,0)[cc]{$\bullet$}}
\put(24,10){\makebox(0,0)[cc]{$\bullet$}}
\put(35,6){\makebox(0,0)[cc]{$\bullet$}}
\put(18,25){\makebox(0,0)[cc]{$f$}}
\put(20,23){\makebox(0,0)[cc]{$a$}}
\put(2,7){\makebox(0,0)[cc]{$b$}}
\put(12,8){\makebox(0,0)[cc]{$c$}}
\put(25,11){\makebox(0,0)[cc]{$d$}}
\put(36,7){\makebox(0,0)[cc]{$v$}}
\put(10,15){\makebox(0,0)[cc]{$g$}}
\put(18,19){\makebox(0,0)[cc]{$h$}}
\put(27,16){\makebox(0,0)[cc]{$i$}}
\put(15,12){\makebox(0,0)[cc]{$j$}}
\put(22.2,14.2){\makebox(0,0)[cc]{$k$}}
\put(20,17){\makebox(0,0)[cc]{$e$}}
\put(0,2){\makebox(0,0)[cc]{$l$}}
\put(3,2){\makebox(0,0)[cc]{$m$}}
\put(6,2){\makebox(0,0)[cc]{$n$}}
\put(9,2){\makebox(0,0)[cc]{$o$}}
\put(13,2){\makebox(0,0)[cc]{$p$}}
\put(17,2){\makebox(0,0)[cc]{$q$}}
\put(19,2){\makebox(0,0)[cc]{$r$}}
\put(24,2){\makebox(0,0)[cc]{$s$}}
\put(32,2){\makebox(0,0)[cc]{$t$}}
\put(38,2){\makebox(0,0)[cc]{$u$}}
\end{picture}
  \caption{Language of trees: a 10-tree.}
  \label{fig:3}
\end{figure}
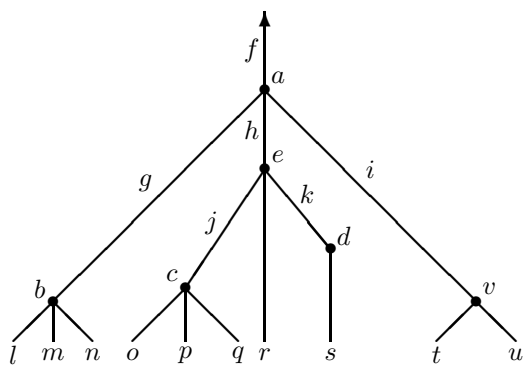
In Figure~\ref{fig:3}, $f$ is the root, $a$ the root vertex,
$l,m,n,o,p,q,r,s,t,u$ the input legs and $g,j,k,i$ the interior edges
of $T$.  Vertices $b,c,d,v$ are the input vertices of $T$.  The vertex $d$
is unary, $v$ binary and the remaining vertices are ternary.  The tree
in Figure~\ref{fig:3} is not reduced because it contains an unary
vertex.  }
\end{example}


\end{document}